\documentclass[letterpaper, 10 pt, conference]{ieeeconf}  

\IEEEoverridecommandlockouts                              
\usepackage{graphicx} 
\usepackage{amsmath} 
\usepackage{amssymb}  
\usepackage{amstext,enumerate}
\usepackage{subfigure}
\usepackage{epstopdf}
\usepackage{url}

\usepackage{latexsym, color,booktabs}

\newtheorem{assum}{Assumption}
\newtheorem{prop}{Proposition}
\newtheorem{thm}{Theorem}

\newtheorem{defn}{Definition}

\title{\LARGE \bf
A secure state estimation algorithm \\ for nonlinear systems under sensor attacks 
}

\author{Michelle S. Chong, Henrik Sandberg, Jo\~{a}o P.~Hespanha 
\thanks{M. Chong is with the Control Systems Technology section at the Department of Mechanical Engineering, Eindhoven University of Technology. 
        {\tt\small m.s.t.chong@tue.nl} }
\thanks{H. Sandberg is with the Division of Decision and Control Systems at KTH Royal Institute of Techonology.
        {\tt\small hsan@kth.se}}
\thanks{J. Hesphana is with the Electrical and Computer Engineering Department at the University of California Santa Barbara.
        {\tt\small hespanha@ece.ucsb.edu}}        
\thanks{This material is based upon work supported by the U.S. Office of Naval Research under the MURI grant No. N00014-16-1-2710.}
}

\begin{document}
\maketitle
\thispagestyle{empty}
\pagestyle{empty}

\begin{abstract}
The state estimation of continuous-time nonlinear systems in which a subset of sensor outputs can be maliciously controlled through injecting a potentially unbounded additive signal is considered in this paper. Analogous to our earlier work for continuous-time linear systems in \cite{chong2015observability}, we term the convergence of the estimates to the true states in the presence of sensor attacks as `observability under $M$ attacks', where $M$ refers to the number of sensors which the attacker has access to. Unlike the linear case, we only provide a sufficient condition such that a nonlinear system is observable under $M$ attacks. The condition requires the existence of asymptotic observers which are robust with respect to the attack signals in an input-to-state stable sense. We show that an algorithm to choose a compatible state estimate from the state estimates generated by the bank of observers achieves asymptotic state reconstruction. We also provide a constructive method for a class of nonlinear systems to design state observers which have the desirable robustness property. The relevance of this study is illustrated on monitoring the safe operation of a power distribution network.  
\end{abstract}

\section{Introduction} \label{sec:intro}

The \textit{cyber} security of dynamical systems have garnered the attention of our community in the past decade, see \cite{sandberg2015cyberphysical} and \cite{chong2019tutorial} for a tutorial overview. This is indeed a timely concern as the increasing (cyber) connectivity between physical systems creates vulnerabilities where malicious cyber attacks can lead to disastrous consequences.

The focus of this work is on the state estimation of nonlinear dynamical systems where the sensors have been compromised. This context has been studied in detail for linear systems in both discrete \cite{fawzi2014secure, shoukry2015event,shoukry2018smt, xie2018secure, lu2019switched, lu2019secure, an2018state} and continuous-time \cite{chong2015observability, an2017secure}. The main thread that underlies achieving state reconstruction is characterizing the number of sensors which are allowed to be attacked, and the resulting algorithm is an optimization problem which is combinatorial in nature. The computational complexity of these algorithms are addressed in various ways including transforming an $l_0$ minimization problem into a convex one \cite{fawzi2014secure}, using gradient descent algorithms \cite{shoukry2015event}, employing Satisfiability Modulo Theory (SMT) solvers to reduce search time \cite{shoukry2018smt} and reducing the number candidates \cite{ an2017secure, an2018state, lu2019secure, lu2019switched}, to name a few. 

Relatively little work has been done for nonlinear systems, where algorithms were proposed for classes of nonlinear systems in discrete-time \cite{shoukry2015secure, hu2017secure, yang2018robust} and in continuous-time \cite{kim2018detection, nateghi2018cyber}. Feedback linearizable systems are considered in \cite{hu2017secure} and differentially flat systems in \cite{shoukry2015secure}, which then enables state estimation using linear techniques. The authors of \cite{yang2018robust} consider Lur'e systems and employs the same framework as in \cite{chong2015observability}. An adaptive observer is designed to estimate both the states and the attack signals for asymptotically stable nonlinear systems in \cite{nateghi2018cyber}. In \cite{kim2018detection}, a uniformly observable nonlinear system is considered and a high gain observer is designed for each measured output. An algorithm which exploits redundancy then collects all the state estimates and provides a state estimate. 

In this paper, we consider a continuous-time nonlinear system with $N$ outputs where each is measured by a potentially compromised sensor. Under the scenario where $M$ out of the $N$ sensors have been maliciously manipulated, we aim to reconstruct the states given that we do not know which of the $M$ sensors have been compromised. If this objective is met, we call such a system \textit{observable under M attacks}, a term coined in our earlier work for linear systems \cite{chong2015observability}.

We first provide a sufficient condition in Section \ref{sec:result} for observability under $M$ attacks. The condition calls for the total number of sensors $N$ to be larger than twice the number of attacked sensors $M$, i.e. $N > 2M$. Moreover, it also requires an observer to be constructed for every combination of $N-2M$ sensor measurements received by the observer, with the crucial property that the observer is robust towards the attack signals. In other words, each observer must have an estimation error system which is input-to-state stable (ISS) \cite{sontag2008input} with respect to the attack signals. These conditions are consistent with the key results in the literature for linear systems \cite{fawzi2014secure, chong2015observability} and a class of nonlinear systems \cite{kim2018detection}.

This gives rise to an algorithm in Section \ref{sec:algo}, which employs the same framework proposed in an earlier work for linear systems in \cite{chong2015observability} by some of the authors of this paper. The algorithm uses a bank of observers designed to satisfy the aforementioned properties and picks the state estimate which satisfies a consistency measure involving a subset of the other state estimates. The chosen state estimate is shown to converge asymptotically to the true state, in the presence of sensor attacks, provided that the system is $M$ attack observable.

In Section \ref{sec:case}, we consider a class of nonlinear systems and provide a systematic method for designing observers which have the desired ISS property with respect to the attack signal. This work is highly relevant in the remote monitoring of the local voltage regulation of each customer who is connected to a power distribution network, which we present in Section \ref{sec:app}. All proofs are provided in the appendix. 
\section{Preliminaries} \label{sec:prelim}
 \begin{itemize}
 	\item Let $\mathbb{R}=(-\infty,\infty)$, $\mathbb{R}_{\geq 0}=[0,\infty)$, $\mathbb{R}_{>0}=(0,\infty)$.
	\item Let the set of complex numbers be denoted by $\mathbb{C}$.
	\item We denote the set of integers $\{i,i+1,i+2,\dots,i+k\}$ as $\mathbb{N}_{[i,i+k]}$.
	\item The number of $k$-element subsets of an $n$-element set is denoted $\binom{n}{k}$. 
 	\item Let $(u,v)$ where $u\in\mathbb{R}^{n_u}$ and $v\in\mathbb{R}^{n_v}$ denote the column vector $(u^T,v^T)^{T}$.
 	\item The cardinality of a set $\mathcal{J}$ is denoted as $\#\left(\mathcal{J}\right)$.
 	\item The identity matrix of dimension $n$ is denoted by $\mathbb{I}_{n}$ and a matrix of dimension $m$ by $n$ with all elements $1$ is denoted by $\mathbf{{1}}_{m\times n}$.
	\item A diagonal matrix with elements $d_i$, $i\in\mathbb{N}_{[1,n]}$ is denoted by $\textrm{diag}(d_1,d_2,\dots,d_n)$.
 	\item Given a symmetric matrix $P$, its maximum (minimum) eigenvalue is denoted by $\lambda_{\max}(P)$ $(\lambda_{\min}(P))$.
	\item The infinity norm of a vector $x \in \mathbb{R}^{n}$, is denoted $|x|:= \underset{i\in\mathbb{N}_{[1,n]}}{\max} \left| x_i \right|$ and for a matrix $A\in\mathbb{R}^{n\times n}$, $|A|:= \underset{i\in\mathbb{N}_{[1,n]}}{\max} \underset{j\in\mathbb{N}_{[1,n]}}{\sum}|a_{ij}|$, where $a_{ij}$ is the row $i$-th and column $j$-th element of matrix $A$.
 	\item A continuous function $\alpha:\mathbb{R}_{\geq 0}\to\mathbb{R}_{\geq 0}$ is a class $\mathcal{K}$ function, if it is strictly increasing and $\alpha(0)=0$; additionally, if $\alpha(r)\to\infty$ as $r\to\infty$, then $\alpha$ is a class $\mathcal{K}_{\infty}$ function. A continuous function $\beta:\mathbb{R}_{\geq0}\times \mathbb{R}_{\geq 0} \to \mathbb{R}_{\geq 0}$ is a class $\mathcal{KL}$ function, if: (i) $\beta(.,s)$ is a class $\mathcal{K}$ function for each $s\geq 0$; (ii) $\beta(r,.)$ is non-increasing and (iii) $\beta(r,s)\to 0$ as $s\to \infty$ for each $r\geq 0$.
 \end{itemize}
\section{Problem statement} \label{sec:prob}
We consider the problem of state observation for a class of nonlinear systems under sensor attacks of the following form
\begin{eqnarray} \label{eq:system}
    \dot{x} & = & f(x,z,w), \, z=(z_1,z_2,\dots,z_N), \nonumber \\
    z_i & = & h_i(x,w), \nonumber \\
    y_i & = & z_i + a_i, \qquad i\in\mathbb{N}_{[1,N]},
\end{eqnarray}
where $x\in\mathbb{R}^{n_x}$ is the state, $y_i\in\mathbb{R}^{n_i}$ is the measured output at sensor $i$, $w\in\mathbb{R}^{n_u}$ is a measured input, $f$ and $h_i$ are locally Lipschitz functions and $a_i\in\mathbb{R}^{n_i}$ is a possibly unbounded attack signal that cannot be measured. 
\begin{assum}
Further assumptions about the attack signals $a_i$ are
\begin{enumerate}[(i)]
    \item Sensors $i\in\mathbb{N}_{[1,N]}$ which are not under attack satisfy $a_i(t)=0$, for all $t\geq 0$.
    \item Given an index set $\mathcal{I}\subseteq \mathbb{N}_{[1,N]}$, the set of non-attacked sensors remain constant, i.e. the attack vector $a=(a_1,a_2,\dots,a_N) \in \mathcal{N}_{I}$, where $\mathcal{N}_{\mathcal{I}}:=\{(a_1,a_2,\dots,a_N): a_i(t)=0, \forall t\geq 0, \forall i \not\in \mathcal{I} \}$.
\end{enumerate} \hfill $\Box$
\end{assum}

In this paper, we derive conditions such that the state $x$ of system \eqref{eq:system} with $N$ outputs can be estimated when $M$ of the sensors have been attacked, which we term \textit{observable under $M$ attacks} and formally define below.

\begin{defn} \label{def:obs_under_attack}
System \eqref{eq:system} is \textbf{observable under $M$ attacks} if for any 
\begin{itemize}
    \item initial conditions $x(0)$, $\bar{x}(0)\in\mathbb{R}^{n_x}$,
    \item measured input $w \in \mathbb{R}^{n_u}$,
    \item index sets $\mathcal{I}_{a}$, $\mathcal{I}_{b} \subset \mathbb{N}_{[1,N]}$ with not more than $M$ elements,
    \item attack vectors $a=(a_1,a_2,\dots,a_{N})\in \mathcal{N}_{\mathcal{I}_a}$, $\bar{a}=(\bar{a}_1,\bar{a}_2,\dots,\bar{a}_{N})\in \mathcal{N}_{\mathcal{I}_b}$,
\end{itemize}
there exists an index set $\mathcal{J} \subset \mathbb{N}_{[1,N]}$ with at least $N-2M$ elements, such that the output trajectories of system \eqref{eq:system} satisfy
\begin{equation}
\begin{aligned} \label{eq:obs_attack_defn}
    \underline{\alpha}_y & \left(\left|x(0)-\bar{x}(0)\right|, t \right)  \\
    & \leq  | y_i\left(t;x(0),w,a_i\right)  -  y_{i} \left(t;\bar{x}(0),w,\bar{a}_i  \right)  |  \\ & \qquad \qquad \leq \bar{\alpha}_y\left(\left|x(0)-\bar{x}(0)\right| \right),
\end{aligned}
\end{equation}
for all $i\in \mathcal{J}$, $t\geq 0$, and $\underline{\alpha}_y$ is a class $\mathcal{KL}$ function and $\bar{\alpha}_y$ is a class $\mathcal{K}$ function. \hfill $\Box$
\end{defn}
We have denoted the output trajectories of system $\eqref{eq:system}$ initialized at $x(0)$ for the input $w$ and attack $a_i$ as $y_i(t;x(0),w,a_i)$, for all $i\in\mathbb{N}_{[1,N]}$.

Definition \ref{def:obs_under_attack} means that when a system \eqref{eq:system} is observable under $M$ attacks, there is at most one initial condition in which system \eqref{eq:system} generates a compatible measured output $y_i$ for any given input signal $w$, for at least $N-2M$ of the measured outputs. This has to be achieved regardless of which of the $M$ sensors have been compromised and the attack signal $a_i$ that has been chosen by the attacker.

\section{A sufficient condition for observabilty under $M$ attacks} \label{sec:result}
We provide a sufficient condition for system \eqref{eq:system} to be observable under $M$ attacks.

\begin{thm} \label{thm:M_obs}
For any integer $M \geq 0$, (ii) implies (i): 
\begin{enumerate}[(i)]
 \item System \eqref{eq:system} is observable under $M$ attacks. 
 \item $N> 2M$ and, for every set $\mathcal{J} \subset \mathbb{N}_{[1,N]}$ with $\#\left(\mathcal{J} \right) \geq N-2M$, there exists a function $\hat{f}:\mathbb{R}^{n_x}\times\mathbb{R}^{\#\left(\mathcal{J} \right)}\times\mathbb{R}^{\#\left(\mathcal{J} \right)} \to \mathbb{R}^{n_x}$ such that the solution to 
 \begin{equation} \label{eq:xhat_gen}
    \dot{\hat{x}}_{\mathcal{J}} = \hat{f}\left( \hat{x}_{\mathcal{J}}, y_{\mathcal{J}}, w \right),  
 \end{equation}
 and the solution to system \eqref{eq:system}, respectively satisfy
\begin{align} \label{eq:obs_iss}
   \left| x(t) - \hat{x}_{\mathcal{J}}(t) \right| \leq  \hat{\beta} &\left(\left| x(0) - \hat{x}_{\mathcal{J}}(0) \right|,t \right) \nonumber \\ & + \hat{\gamma}\left( \underset{s\in[0,t)}{\sup} \left|a_{\mathcal{J}}(s) \right| \right), 
\end{align}
for all $t\geq 0$ and initial conditions $x(0)$, $\hat{x}_{\mathcal{J}}(0) \in \mathbb{R}^{n_x}$, where $\hat{\beta}$ is a $\mathcal{KL}$ function, $\hat{\gamma}$ is a $\mathcal{K}_{\infty}$ function, and $a_{\mathcal{J}}$ denotes a stacked vector of $a_i$ indexed by $i\in\mathcal{J}$. 
\end{enumerate} \hfill $\Box$
\end{thm}

Theorem \ref{thm:M_obs} specifies that the number of available sensors $N$ has to be strictly more than twice the number of compromised sensors $M$. This is consistent with the results for linear systems in \cite{chong2015observability} for continuous-time systems and \cite{fawzi2014secure} for discrete-time systems, as well as in \cite{kim2018detection} and \cite{yang2018robust} for classes of nonlinear systems in continuous and discrete-time, respectively. 

Further, condition (ii) means that the estimation error $e_{\mathcal{J}}:=x-\hat{x}_{\mathcal{J}}$ system constructed out of system \eqref{eq:system} and \eqref{eq:xhat_gen} is input-to-state stable (ISS) \cite{sontag2008input} with respect to the attack vector $a_{\mathcal{J}}$. This property can be fulfilled with Luenberger observers in the case of linear systems (see \cite[Section III.B]{chong2015observability}), and with high gain observers \cite{bornard1991high} or circle criterion observers \cite{arcak2001nonlinear} for classes of nonlinear systems. We will provide a constructive example in our case study in Section \ref{sec:case}.


\section{Algorithm} \label{sec:algo}
Using Theorem \ref{thm:M_obs}, we formulate the following estimation algorithm to estimate the states of system \eqref{eq:system} when $M$ out of $N$ of its sensors have been compromised. Our algorithm follows the idea presented in \cite{chong2015observability}, where the results were derived for linear dynamical systems. 

The crux of the algorithm lies in the fact that for each combination of $N-M$ outputs (note that this is greater than $N-2M$ outputs, which satisfies condition (ii) in Theorem~\ref{thm:M_obs}), there is one observer which receives attack-free sensor outputs and hence provides state estimates that converges to the true state. Further, for each set of these $N-M$ outputs, there is at least one subset consisting of $N-2M$ outputs which is attack-free. Thus, the observer which receives the attack-free subset of $N-2M$ outputs will provide a state estimate which converge to the true one. Therefore, in the algorithm presented in this section, we employ two banks of observers: one bank of observers employing $N-M$ outputs, and the other employing $N-2M$ outputs. 

Suppose that at most $M$ out of $N$ of system \eqref{eq:system}'s outputs can be compromised and condition (ii) of Theorem~\ref{eq:system} holds. Then, for every set $\mathcal{S} \subset \mathbb{N}_{[1,N]}$ of $N-M$ elements, an observer which employs $N-M$ outputs from system \eqref{eq:system} is constructed as follows
\begin{equation} \label{eq:obs_gen_algo}
    \dot{\hat{x}}_{\mathcal{S}} = \hat{f}(\hat{x}_{\mathcal{S}},y_{\mathcal{S}},w),
\end{equation}
which has an estimation error system that is ISS with respect to the attack vector $a_{\mathcal{S}}$ as stated in \eqref{eq:obs_iss} of Theorem \ref{thm:M_obs}. This forms the first bank of $\binom{N}{N-M}$ observers. We define the consistency measure $\pi_{\mathcal{S}}$ to be the worst case deviation between the estimate $\hat{x}_{\mathcal{S}}$ given by \eqref{eq:obs_gen_algo} and the estimate $\hat{x}_{\mathcal{P}}$ generated in the same manner as \eqref{eq:obs_gen_algo} for $\mathcal{P}\subset \mathcal{S}$ with $N-2M$ elements, which is
\begin{equation} \label{eq:cons_measure}
    \pi_{\mathcal{S}}(t) = \underset{\mathcal{P}\subset\mathcal{S}: \#(\mathcal{P})=N-2M}{\max} \left| \hat{x}_{\mathcal{S}}(t) - \hat{x}_{\mathcal{P}}(t) \right|.
\end{equation}

For the set $\mathcal{S}$ in which all the attack vectors are zero, i.e. $a_i(t)=0$, for all $i\in\mathcal{S}$ and $t\geq 0$, all the state estimates $\hat{x}_{\mathcal{S}}$ and $\hat{x}_{\mathcal{P}}$ will be consistent and this motivates the choice of the state estimate $\hat{x}$ produced by the algorithm as follows
\begin{align} \label{eq:algo_choice}
    \hat{x}(t) & = \hat{x}_{\sigma(t)}(t), \nonumber \\
    \sigma(t) & = \underset{\mathcal{S}\subset\mathbb{N}_{[1,N]}:\#(\mathcal{S})=N-M}{\arg\min} \pi_{\mathcal{S}}(t).
\end{align}

\begin{figure}
    \centering
    \includegraphics{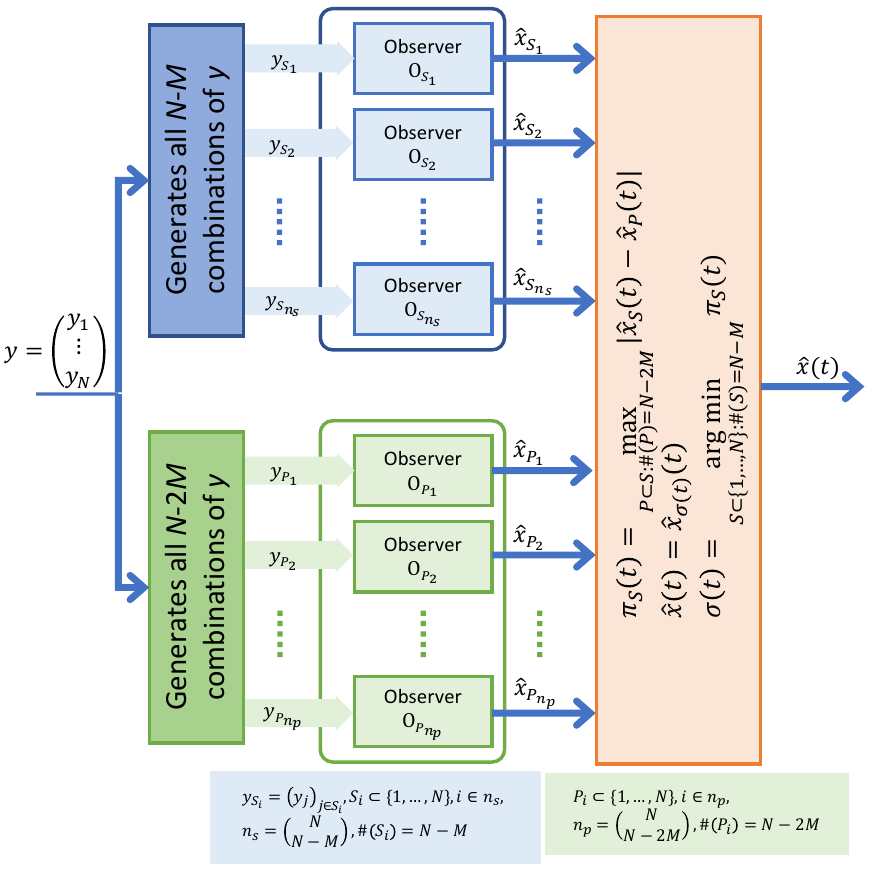}
    \caption{Infrastructure of a secure state estimation algorithm for nonlinear systems \eqref{eq:system}.}
    \label{fig:secure_obs}
\end{figure}

We summarize the algorithm \eqref{eq:obs_gen_algo}, \eqref{eq:cons_measure}, \eqref{eq:algo_choice} in Figure \ref{fig:secure_obs} and provide the following state estimation convergence guarantees. 

\begin{thm} \label{thm:algo}
Consider system \eqref{eq:system} with $N$-outputs of which at most $M$ is compromised, i.e. the attack vector $a$ belongs to $\mathcal{N}_{\mathcal{I}}$, for some set $\mathcal{I}\subset\mathbb{N}_{[1,N]}$ where $\#(\mathcal{I})\leq M$. Assuming that (ii) of Theorem \ref{thm:M_obs} holds, then there exists a class $\mathcal{KL}$ function $\tilde{\beta}$ such that the solution to system \eqref{eq:system} and the secure state estimation algorithm \eqref{eq:obs_gen_algo}, \eqref{eq:cons_measure}, \eqref{eq:algo_choice} satisfy
\begin{equation} \label{eq:tildex_algo}
    \left|x(t)-\hat{x}(t)\right| \leq \tilde{\beta}\left(\left|x(0)-\hat{x}(0)\right|,t\right), \; \forall t\geq 0,
\end{equation}
for any initial conditions $x(0)$, $\hat{x}_{\mathcal{S}}(0)$, $\hat{x}_{\mathcal{P}}(0) \in \mathbb{R}^{n_x}$. \hfill $\Box$
\end{thm}

\section{Case study: A class of nonlinear systems} \label{sec:case}

We consider a specific form of system \eqref{eq:system} as follows:
\begin{eqnarray} \label{eq:system_spec}
    \dot{x} & = & Ax + \phi(z), \; \phi(z)=\left(\phi_1(z_1),\phi_2(z_2),\dots,\phi_N(z_N) \right), \nonumber \\
    z_i & = & H_i x + w_i, \; i\in\mathbb{N}_{[1,N]}, \nonumber \\
    y_i & = & z_i + a_i,
\end{eqnarray}
where the nonlinearities $\phi_i:\mathbb{R}^{n_i} \to \mathbb{R}$ are slope-restricted, i.e.
\begin{assum} \label{assum:sector}
For $i\in\mathbb{N}_{[1,N]}$, the nonlinearity $\phi_i$ satisfies
\begin{equation} \label{eq:sector_bounded}
    \underline{d}_{i} \leq \frac{\phi_i(\xi)-\phi_i(\psi)}{\xi-\psi} \leq \bar{d}_i, \; \forall \xi, \psi \in\mathbb{R},\, \xi \neq \psi.
\end{equation} \hfill $\Box$
\end{assum}

For system \eqref{eq:system_spec} with $N$ outputs, of which $M$ can be compromised, we show that (ii) of Theorem \ref{thm:M_obs} is satisfied by designing each observer \eqref{eq:xhat_gen} in the following manner for every set $\mathcal{J}\subset\mathbb{N}_{[1,N]}$.
\begin{eqnarray} \label{eq:xhat_J}
    \dot{\hat{x}}_{\mathcal{J}} & = & A \hat{x}_{\mathcal{J}} + \phi(\xi_{\mathcal{J}}) + \bar{o}(L_{\mathcal{J}},y_{\mathcal{J}},\hat{x}_{\mathcal{J}},w_{\mathcal{J}}), \nonumber \\
    \xi_{\mathcal{J}} & = & H\hat{x}_{\mathcal{J}}+w+\bar{o}(K_{\mathcal{J}},y_{\mathcal{J}},\hat{x}_{\mathcal{J}},w_{\mathcal{J}}),
\end{eqnarray}
where $\bar{o}(K_{\mathcal{J}},y_{\mathcal{J}},\hat{x}_{\mathcal{J}},w_{\mathcal{J}}) = K_{\mathcal{J}} \left(y_{\mathcal{J}} - \left( H_{\mathcal{J}}\hat{x}_{\mathcal{J}} + w_{\mathcal{J}} \right) \right)$ is an output injection term employing outputs $y_{\mathcal{J}}$, known inputs $w_{\mathcal{J}}$ and an observer matrix $K_{\mathcal{J}}$ to be designed. Note that the first two terms in $\xi_{\mathcal{J}}$ use the full $H$ from system \eqref{eq:system_spec} and all the known inputs $w$, respectively. 

\begin{prop} \label{prop:spec}
    Consider system \eqref{eq:system_spec} under Assumption \ref{assum:sector}. Suppose $N> 2M$ and, for every set $\mathcal{J} \subset \mathbb{N}_{[1,N]}$ with $\#\left(\mathcal{J} \right) \geq N-2M$, there exist a matrix $P_{\mathcal{J}}=P_{\mathcal{J}}^{T}>0$, scalars $\nu_{\mathcal{J}} \geq 0$, $\mu_{\mathcal{J}} \geq 0$ and observer matrices $K_{\mathcal{J}}$ and $L_{\mathcal{J}}$ such that the following holds
    \begin{equation} \label{eq:cc_lmi}
        \left[\begin{array}{ccc} \mathcal{A}\left(P_{\mathcal{J}}, P_{\mathcal{J}}L_{\mathcal{J}}, \nu_{\mathcal{J}} \right) & \mathcal{B}\left(P_{\mathcal{J}}, K_{\mathcal{J}}  \right) & -P_{\mathcal{J}} \\
        \mathcal{B}\left(P_{\mathcal{J}}, K_{\mathcal{J}}  \right)^{T} & \mathcal{D}(\bar{d}) & 0 \\
        -P_{\mathcal{J}} & 0 & -\mu_{\mathcal{J}} \mathbb{I}_{n_{\mathcal{J}}} \end{array} \right] \leq 0,
    \end{equation}
    where 
    \begin{itemize}
        \item $\mathcal{A}\left(P_{\mathcal{J}}, P_{\mathcal{J}}L_{\mathcal{J}}, \nu_{\mathcal{J}} \right):=$ \\ $ P_{\mathcal{J}}\left(A - L_{\mathcal{J}} H_{\mathcal{J}} \right) + \left(A - L_{\mathcal{J}} H_{\mathcal{J}} \right)^{T} P_{\mathcal{J}} + \nu_{\mathcal{J}} \mathbb{I}_{n_x}$, 
        \item $\mathcal{B}\left(P_{\mathcal{J}}, K_{\mathcal{J}}  \right):=P_{\mathcal{J}}+\left(H-K_{\mathcal{J}}H_{\mathcal{J}}\right)^{T}$,
        \item $\mathcal{D}(\bar{d}):=-2 \, \textrm{diag}\left(\bar{d}_{1}^{-1},\bar{d}_{2}^{-1},\dots, \bar{d}_{N}^{-1}\right)$,
        \item $n_{\mathcal{J}}:=\underset{i\in\mathcal{J}}{\sum} n_i$.
    \end{itemize}
    Then, (ii) of Theorem \ref{thm:M_obs} holds. \hfill $\Box$
\end{prop}

Inequality \eqref{eq:cc_lmi} is a linear matrix inequality (LMI) in $P_{\mathcal{J}}$, $P_{\mathcal{J}}L_{\mathcal{J}}$, $
\nu_{\mathcal{J}}$, $K_{\mathcal{J}}$ and $\mu_{\mathcal{J}}$, which can be solved efficiently using computational tools. The design we have used here was first introduced as the circle criterion observer in \cite{arcak2001nonlinear}, which can be tuned to attenuate measurement noise and input disturbances according to the design in \cite{chong2012robust}. Here, we have adapted the design such that the observer \eqref{eq:xhat_J} is robust with respect to the attack vector $a_{\mathcal{J}}$ in the sense of \eqref{eq:obs_iss}.

\section{Application: Secure monitoring for the voltage regulation of a power distribution network} \label{sec:app}

\begin{figure}[h!]
    \centering
    \includegraphics{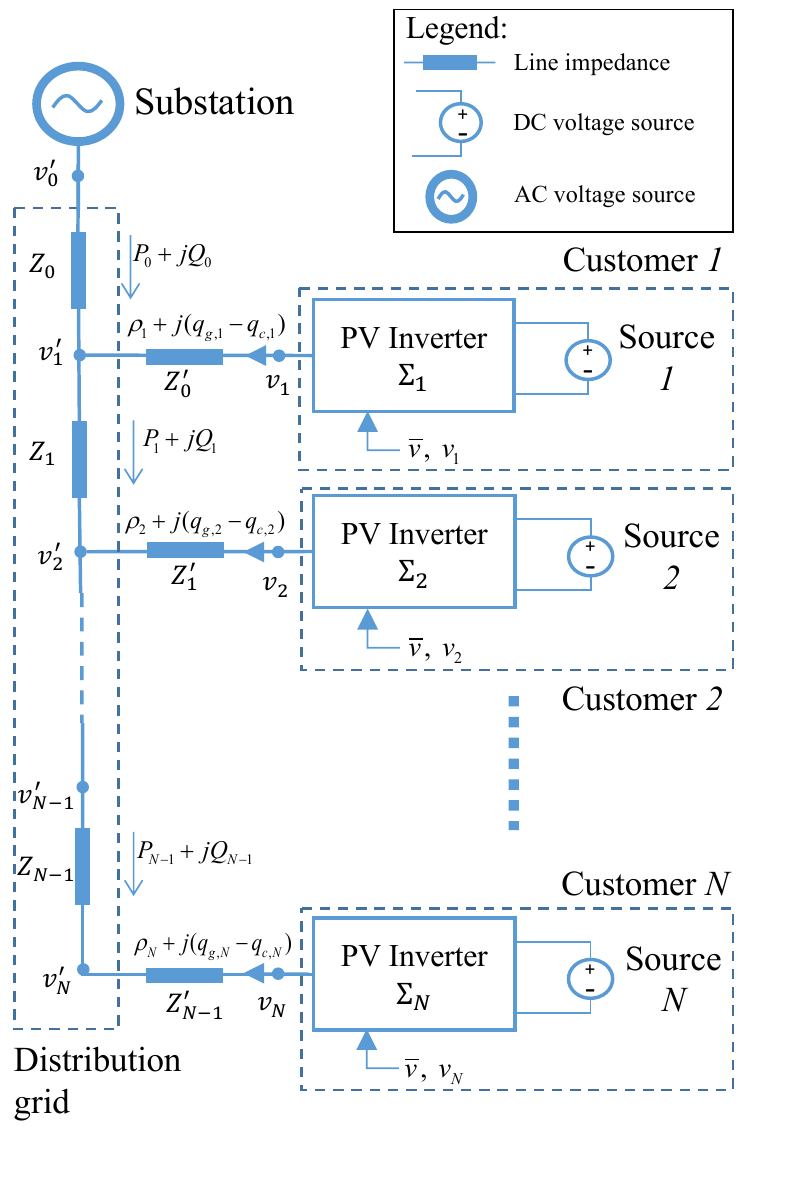} \vspace{-1em}
    \caption{Infrastructure of a low-voltage grid. Voltage regulation is achieved via local controllers $\Sigma_i$, for $i\in\mathbb{N}_{[1,N]}$, with the monitoring center (situated at the substation) receiving potentially corrupted measurements of $v_i$.}
    \label{fig:lv_grid}
\end{figure}

A typical low voltage power distribution network (shown in Figure \ref{fig:lv_grid}) would consist of $N$ customers feeding into the distribution network in a line configuration, with the smart secondary substation at the head of the line. The substation functions as a monitoring center, sending the desired set-point voltage $\bar{v}$ to each local controller $\Sigma_i$, such that the voltages received by each customer $v_i$ is regulated to operate in a safe operating range, i.e. for a given $\delta>0$,
\begin{equation} 
    \bar{v}-\delta \leq v_i(t) \leq \bar{v} + \delta, \qquad \forall t\geq 0. \label{eq:voltage_safe}
\end{equation}

For each customer $i\in\mathbb{N}_{[1,N]}$, the received voltage level is $v_i$ and the voltage level at the point of connection with the distribution line is $v'_i$, with a corresponding line impedance $Z'_i=R'_{i} + jX'_{i}$ in between customer $i$ and the connection point on the distribution line, where $R'_{i}\in\mathbb{R}_{\geq 0}$ is the resistance and $X'_{i}\in\mathbb{R}_{\geq 0}$ is the reactance. In between each connection point, the corresponding line impedance is $Z_{i}=R_{i} + jX_{i}$, where $R_{i}\in\mathbb{R}_{\geq 0}$ is the resistance and $X_{i}\in\mathbb{R}_{\geq 0}$ is the reactance. Each customer has a load which can consume reactive $q_{c,i}$ and active powers $\rho_{c,i}$, independently of the generated reactive $q_{g,i}$ and active powers $\rho_{g,i}$.

In \cite{chong2019local}, a class of sector-bounded droop controllers $\Sigma_i$ which uses local measurements $v_i$ were shown to regulate the voltages $v_i$ such that the safety constraint \eqref{eq:voltage_safe} is satisfied. This is achieved via appropriate injection of reactive power $q_{g,i}$ by each local controller $\Sigma_i$ to regulate the flow of active $P_i$ and reactive $Q_i$ powers, under the assumption that the net injected active power $\rho_i$ and the reactive power $q_{c,i}$ consumed by customer $i$ is bounded and the bounds are known. 

We are now concerned with the \textit{security} problem where the measurements of the voltages $v_i$ received at the monitoring center situated at the substation has been maliciously corrupted. We model this measurement corruption with an additive attack signal $\alpha_i:\mathbb{R}_{\geq 0} \to \mathbb{R}$ which is potentially unbounded, as follows
\begin{equation} \label{eq:corrupt_voltage}
    \hat{y}_i = v_i + \alpha_i.
\end{equation}
This is a major issue as the presence of the attack signal $\alpha_i$ would mislead the monitoring center into thinking that the safety constraint \eqref{eq:voltage_safe} has been violated and thus a false alarm is raised and possibly triggering unwarranted operator actions. 

Our solution is to employ the results in the previous sections to estimate the voltages $v_i$, given that $M$ out of $N$ of the measurements $y_i$ are maliciously manipulated. To this end, we model the power flow in the distribution grid as done in \cite{chong2019local} using the linearized DistFlow model \cite{baran1989network} and assume that the droop controllers $\Sigma_i$ have been designed using the methodology presented in \cite{chong2019local}. The relationship between the power flow and voltages between key nodes is
\begin{equation} \label{eq:line_rel}
	\begin{array}{lll}
			P_{i+1} & = & P_{i} + \rho_{i+1}, \\
			Q_{i+1} & = & Q_{i} + q_{i+1}, \\
			{v'_{i+1}}^2 & = & {v'_{i}}^2 - 2\beta_{i}(P_i,Q_i), \\
			{v'_i}^2 & = & {{v}_{i}}^2 - 2 {\beta}^{'}_{i-1}(\rho_i,q_i),
	\end{array}
\end{equation}
where $P_i$ and $Q_i$ are the respective total active and reactive powers flowing from customer $i$ to customer $i+1$; $\rho_i:=\rho_{g,i}-\rho_{c,i}$ and $q_i:=q_{g,i} - q_{c,i}$ are the net injection of the respective active and reactive power into the distribution line from customer $i$; $\beta_{i}(r,s):=R_i r + X_i s$ and $\beta'_{i}(r,s):=R'_i r + X'_i s$ with $\beta'_{-1}(r,s)=0$ for all $r,s\in\mathbb{R}$. 

Each local controller $\Sigma_i$ actuated by the inverter can generate reactive power $q_{g,i}$ as follows
\begin{equation} \label{eq:droop}
		\dot{q}_{g,i} = -\frac{1}{\tau_i} q_{g,i} + \frac{1}{\tau_i} K_{i}(\bar{v}^2-v_{i}^2), 
\end{equation} 
where $\tau_i\in\mathbb{R}_{> 0}$ is the time-constant of the inverter's response, $\bar{v}\in\mathbb{R}$ is the reference voltage communicated to each customer $i$ and the droop function $K_{i}(w)$ is a static mapping from the difference of the squared voltages $w$ to the set-point for the reactive power. We choose the droop function $K_i(w)$ to be a piecewise saturation function considered in \cite{andren2015stability} which takes the following form:
	\begin{equation}\label{eq:droop_K}
			K_i(w) := \left\{\begin{array}{ll} -\bar{Q}_{i}, & w \leq w_{\min,i}, \\  
			  -\left(1-  \frac{w-w_{\min,i}}{w_{m,i}-w_{\min,i}}\right) \bar{Q}_{i}  & w\in (w_{\min,i},w_{m,i}], \\
			 0, & w\in (w_{m,i},w_{n,i}], \\
		  \left(\frac{w-w_{n}}{w_{\max,i}-w_{n,i}}\right) \bar{Q}_{i}  &  w\in (w_{n,i},w_{\max,i}], \\
			 \bar{Q}_{i}, & w > w_{\max,i},  \end{array} \right. 
	\end{equation}
where $w_{\min,i} \leq w_{m,i} \leq 0 \leq w_{n,i} \leq w_{\max,i}$ are design parameters, $\bar{Q}_{i}\in\mathbb{R}_{\geq 0}$ is the saturation limit of the $i$-th inverter satisfying $\bar{Q}_i=\sqrt{\bar{s}_i^2-\rho_{g,i}^2}$, where $\bar{s}_i\in\mathbb{R}$ is the maximum apparent power of the $i$-th inverter. The design parameters $w_{\min,i}$, $w_{\max,i}$, $w_{m,i}$, $w_{n,i}$ are chosen such that
\begin{equation} \label{eq:di}
		d_i:=\min \left\{\begin{array}{cc} \frac{\bar{Q}_{i}}{w_{\max,i} - w_{n,i}}, \frac{\bar{Q}_{i}}{w_{m,i} - w_{\min,i}} \end{array} \right\}
\end{equation}
satisfies \cite[Theorem 6]{chong2019local} such that the safety constraint \eqref{eq:voltage_safe} is met. We employ the same change in state coordinates as done in \cite{chong2019local} such that the distribution model \eqref{eq:line_rel}, controllers \eqref{eq:droop}, and measurements received at the monitoring center \eqref{eq:corrupt_voltage} can be written in the form of \eqref{eq:system_spec} by choosing
\begin{itemize}
    \item the state $x=\left(q_{g,1},q_{g,2},\dots,q_{g,N}\right)$,
    \item $z_i := {\bar{v}}^2-v_i^2$,
    \item the known input $w_i = \phi_i(\rho,q_c)+{\bar{v}}^{2}-{v'}_{0}^2$, where $\phi_i(\rho, q_c):=\underset{j\in\mathbb{N}_{[0,i-1]}}{\sum} \psi_{j}(\rho,q_c)+\!\!\! \underset{j\in\mathbb{N}_{[0,i-2]}}{\sum} 2{\beta}^{'}_{j}(\rho_{j+1},q_{c,j+1})$ and \begin{equation*}
			\begin{aligned}
				\psi_{j}(\rho,q_c)&:=  2 X_j  \underset{k\in\mathbb{N}_{[j+1,N]}}{\sum} q_{c,k} -2 R_j   \underset{k\in\mathbb{N}_{[j+1,N]}}{\sum} \rho_k \\ &- 2 {\beta'}_{j}(\rho_{j+1},q_{c,j+1}),
			\end{aligned}
		\end{equation*}
		where $\beta'_{j}$ is from \eqref{eq:line_rel},
    \item the attack signal $a_i = 2v_i \alpha_i - \alpha_i^2$, where $\alpha_i$ comes from \eqref{eq:corrupt_voltage},
    \item $H_i$ to be the rows of the matrix 
    \begin{equation*}
    \begin{aligned}
				H &= -2\left(\begin{array}{ccccc} X_0 &  X_0 & \dots  & X_0 \\
			 								  \star & \left(X_0+X_1\right) & \dots  & \left(X_0+X_1 \right) \\
											  \vdots & \ddots & \ddots & \vdots \\
											  \star & \dots & \star & \underset{i\in\mathbb{N}_{[0,N-1]}}{\sum} X_i \end{array}\right)\\ 
											  &-2  \textrm{diag}\left(X'_0,\dots,X'_{N-1} \right),
			\end{aligned}
		\end{equation*}				
		where $\star$ denotes a block component of a symmetric matrix,
    \item $A=\textrm{diag}\left(-1/\tau_{1},-1/\tau_{2},\dots,-1/\tau_{N}\right)$,
    \item $\phi_i(z_i) = \tau_i^{-1} K_i(z_i)$, which satisfies Assumption \ref{assum:sector} with $\underline{d}_{i}=0$ and $\bar{d}_{i}= d_i/\tau_i $, where $d_i$ is defined in \eqref{eq:di}.
\end{itemize} 

To recapitulate, given that the monitoring center sits remotely at the substation, the objective is to estimate the voltages $v_i$, given that the monitoring center only has access to the measurements $y_i$, $i\in\mathbb{N}_{[1,N]}$, where $M$ out of $N$ of these measurements may be corrupted. The main idea is to first estimate the states $x$ of all the controllers $\Sigma_i$, $i\in\mathbb{N}_{[1,N]}$,  then estimate the voltages $v_i$ via
\begin{equation} \label{eq:vhat}
    \hat{v}_i(t)^2 = - H_i \hat{x}(t) - \psi_i(\rho,q_c) + {v'}^{2}_{0},   
\end{equation}
where $\hat{x}$ is the state estimate provided by the secure estimation algorithm described in Section \ref{sec:algo};  $\psi$ and ${v'}^{2}_{0}$ are known. We have kept the squared form of the voltage $\hat{v}_i^2$ due to the distribution model \eqref{eq:line_rel} used. We provide the following guarantee.

\begin{prop} \label{prop:voltage}
    Consider the distribution model \eqref{eq:line_rel} and controllers \eqref{eq:droop} with measurements \eqref{eq:corrupt_voltage} where $a_{\mathcal{I}}$ belongs to $\mathcal{N}_{\mathcal{I}}$ for some unknown set $\mathcal{I}\subset \mathbb{N}_{[1,N]}$ with at most $M$ elements. Suppose $N>2M$ and for every $\mathcal{J}\subset \mathbb{N}_{[1,N]}$ with $\#\left(\mathcal{J}\right)\geq N-2M$, there exist  a matrix $P_{\mathcal{J}}=P_{\mathcal{J}}^{T}>0$, scalars $\nu \geq 0$, $\mu_{a} \geq 0$ and observer matrices $K_{\mathcal{J}}$ and $L_{\mathcal{J}}$ such that \eqref{eq:cc_lmi} holds. Then, using the secure state estimation algorithm \eqref{eq:obs_gen_algo}, \eqref{eq:cons_measure}, \eqref{eq:algo_choice}, the estimated squared voltages $\hat{v}_i^2$ computed according to \eqref{eq:vhat} converges to to the true squared voltages $v_{i}^{2}$ as follows 
    \begin{equation} \label{eq:v2_vhat2}
        \left| v_{i}(t)^2 - \hat{v}_{i}(t)^2 \right| \leq \beta_{v}\left(\left|x(0)-\hat{x}(0)\right|,t\right),
    \end{equation}
    for all $t\geq 0$, $i\in\mathbb{N}_{[1,N]}$,   initial conditions $q_{g,i}(0)\in\mathbb{R}$, where $x(0)=\left(q_{g,1}(0),q_{g,2}(0),\dots,q_{g,N}(0)\right)$, $\hat{x}(0)=\left(\hat{q}_{g,1}(0),\hat{q}_{g,2}(0),\dots,\hat{q}_{g,N}(0)\right)$ and $\beta_{v}$ is a class $\mathcal{KL}$ function. 
    
    \hfill $\Box$
\end{prop}
\section{Conclusions and future work} \label{sec:conc}
We introduced a new definition of observability for nonlinear systems in which a subset of the outputs can be manipulated maliciously. A sufficient condition such that asymptotic state reconstruction can be achieved in the presence of sensor attacks is provided, which requires building a bank of observers with an ISS property. A secure state estimation algorithm is proposed which shows that the framework used for linear continuous time systems in our earlier work \cite{chong2015observability} can be used for nonlinear systems as well. A systematic method for designing the observers is proposed for a class of nonlinear systems and we showed the relevance of this work in the monitoring of a power distribution network. Future work includes reducing the computational resources needed in terms of time and number of observers required.    


\appendix
\subsection{Proof of Theorem \ref{thm:M_obs}}
\textbf{(ii) implies (i)}: Suppose to the contrary that (ii) is true, but (i) is false, i.e. there exist initial conditions $x(0)$, $\bar{x}(0)\in\mathbb{R}^{n_x}$, sets $\mathcal{I}_{a}$, $\mathcal{I}_{b} \subset \mathbb{N}_{[1,N]}$ where $\#\left(\mathcal{I}_{a}\right) \leq M$ and $\# \left(\mathcal{I}_{b}\right) \leq M$, input $w$ and attack vectors $a
\in \mathcal{N}_{\mathcal{I}_{a}}$, $\bar{a}
\in \mathcal{N}_{\mathcal{I}_{b}}$ and an index set $\mathcal{J}\subset \mathbb{N}_{[1,N]}$ with $\# \left( \mathcal{J} \right)\geq N-2M$ such that for all $i\in \mathcal{J}$
\begin{align} \label{eq:false_i_1}
    | y_i\left(t;x(0),w,a_i \right) & - y_i\left(t;\bar{x}(0),w,\bar{a}_{i}\right) | \nonumber \\  & > \bar{\alpha}_y\left(\left|x(0)-\bar{x}(0) \right|\right), \qquad \forall t\geq 0,
\end{align}
and/or 
\begin{align} \label{eq:false_i_2}
    | y_i\left(t;x(0),w,a_i \right) & - y_i\left(t;\bar{x}(0),w,\bar{a}_{i}\right) | \nonumber \\  & < \underline{\alpha}_y\left(\left|x(0)-\bar{x}(0) \right|, t\right), \qquad \forall t\geq 0,
\end{align}
where $\underline{\alpha}_y$ is a $\mathcal{KL}$ function and $\bar{\alpha}$ is a $\mathcal{K}_{\infty}$ function.

First, note that the solution to \eqref{eq:system} satisfies the following for $i\in \mathbb{N}_{[1,N]}$
\begin{align} \label{eq:yi_inter}
    | y_i(t;x(0),&w,a_i )  - y_i\left(t;\bar{x}(0),w,\bar{a}_{i}\right) | \nonumber \\
    & \leq l_i  \left| x(t;x(0),w,a_i) - x(t;\bar{x}(0),w,\bar{a}_i) \right| ,
\end{align}
where $l_i\geq 0$ is the Lipschitz constant of the output function $h_i$ from system \eqref{eq:system}.

Consider the index set $\mathcal{J}:=\mathbb{N}_{[1,N]}/\{\mathcal{I}_{a} \cup \mathcal{I}_{b} \}$, where $a_i=0$ and $\bar{a}_{i}=0$, for all $i\in\mathcal{J}$. 
We can view $x(t;\bar{x}(0),w,\bar{a}_i)$ as a solution to \eqref{eq:xhat_gen} with $\hat{f}\left(\hat{x}_{\mathcal{J}},y_{\mathcal{J}}|_{a_{\mathcal{J}}=0}, w \right) = f(x,z,w)$, which is a copy of system \eqref{eq:system}'s dynamics, initialized at $\bar{x}(0)$. Using (ii), we have from \eqref{eq:yi_inter} that for all $i\in\mathcal{J}$
\begin{align} \label{eq:yi_inter_final}
    | y_i\left(t;x(0),w,a_i \right) & - y_i\left(t;\bar{x}(0),w,\bar{a}_{i}\right) | \nonumber \\  & \leq  \alpha_i \left(\left|x(0)-\bar{x}(0) \right| \right), \qquad \forall t\geq 0,
\end{align}
where $\alpha_i(r) := l_i \hat{\beta}(r,0)$ is a $ \mathcal{K}_{\infty}$ function and $\hat{\beta}$ is a $ \mathcal{KL}$ from (ii). According to \eqref{eq:yi_inter_final}, $| y_i\left(t;x(0),w,a_i \right) - y_i\left(t;\bar{x}(0),w,\bar{a}_{i}\right) | = 0$ when $|x(0)-\bar{x}(0)|=0$, which contradicts \eqref{eq:false_i_1} and \eqref{eq:false_i_2}. Therefore, (ii) implies (i).   

\subsection{Proof of Theorem \ref{thm:algo}}
Let $\bar{\mathcal{I}}:=\mathbb{N}_{[1,N]}\backslash \mathcal{I}$ and consider the consistency measure for the set $\bar{\mathcal{I}}$. 
\begin{align}
    \pi_{\bar{\mathcal{I}}}(t) & =  \underset{\mathcal{P}\subset \bar{\mathcal{I}}:\#(\mathcal{P})=N-2M}{\max} \left| \hat{x}_{\bar{\mathcal{I}}}(t) - \hat{x}_{\mathcal{P}}(t)  \right| \nonumber \\
    & =  \underset{\mathcal{P}\subset     \bar{\mathcal{I}}:\#(\mathcal{P})=N-2M}{\max} \left| \hat{x}_{\bar{\mathcal{I}}}(t) - x(t) + x(t) - \hat{x}_{\mathcal{P}}(t)  \right| \nonumber \\
    & \leq   \left| \hat{x}_{\bar{\mathcal{I}}}(t) - x(t) \right| + \underset{\mathcal{P}\subset     \bar{\mathcal{I}}:\#(\mathcal{P})=N-2M}{\max} \left| x(t) - \hat{x}_{\mathcal{P}}(t)  \right|. 
    \label{eq:pi_Ibar}
\end{align}

Note that $a_i(t)=0$, for all $t\geq 0$ and $i\in\bar{\mathcal{I}}$. Further, for every $\mathcal{P} \subset \bar{\mathcal{I}}$, where $\#(\mathcal{P})=N-2M$, we also have that $a_i(t)=0$, for all $t\geq 0$ and $i\in{\mathcal{P}}$. Therefore, by assumption that \eqref{eq:obs_iss} of Theorem \ref{thm:M_obs} holds, we have that for any $\mathcal{X} \subset \mathbb{N}_{[1,N]}$ (either $\mathcal{X}=\bar{\mathcal{I}}$ or $\mathcal{X}={\mathcal{P}}$) 
\begin{equation}
    \begin{array}{lll}
        \left| \hat{x}_{\mathcal{X}}(t)-x(t) \right| & \leq & \beta_{\mathcal{X}}\left(\left| \hat{x}_{\mathcal{X}}(0)-x(0) \right|,t \right), \; \forall t\geq 0,
    \end{array} \label{eq:tildex_X}
\end{equation}
for all $\hat{x}_{\mathcal{X}}(0)$, $x(0)\in\mathbb{R}^{n_x}$, where $\beta_{\mathcal{X}}$ is a class $\mathcal{KL}$ function.

We also observe that for every set $\mathcal{S}$ with $\#(\mathcal{S})=N-M$, we have at least one set $\bar{\mathcal{P}}\subset \mathcal{S}$ with $\#(\bar{\mathcal{P}}) = N-2M$ where $a_i(t)=0$, for all $t\geq 0$ and $i\in\bar{\mathcal{P}}$. Hence, assuming that \eqref{eq:obs_iss} of Theorem \ref{thm:M_obs} is satisfied, we have that 
\begin{equation}
    \begin{array}{lll}
        \left| \hat{x}_{\bar{\mathcal{P}}}(t)-x(t) \right| & \leq & \beta_{\bar{\mathcal{P}}}\left(\left| \hat{x}_{\bar{\mathcal{P}}}(0)-x(0) \right|,t \right), \; \forall t\geq 0,
    \end{array} \label{eq:tildex_Pbar}
\end{equation}
for all $\hat{x}_{\bar{\mathcal{P}}}(0)$, $x(0)\in\mathbb{R}^{n_x}$, where $\beta_{\bar{\mathcal{P}}}$ is a class $\mathcal{KL}$ function. 

Moreover, using the fact that
\begin{equation}
    \begin{array}{lll}
    \pi_{\sigma(t)}(t) & = & \underset{\mathcal{P}\subset \sigma:\#(\mathcal{P})=N-2M}{\max} \left| \hat{x}_{\sigma(t)}(t) - \hat{x}_{\mathcal{P}}(t)  \right| \\
    & \geq & \left| \hat{x}_{\sigma(t)} - \hat{x}_{\bar{\mathcal{P}}} (t) \right|,
    \end{array} \label{eq:pi_sigma}
\end{equation}
we obtain
\begin{equation}
    \begin{array}{lll}
        \left|x(t)-\hat{x}_{\sigma(t)}(t)\right| & = & \left|x(t)-\hat{x}_{\bar{\mathcal{P}}} + \hat{x}_{\bar{\mathcal{P}}} - \hat{x}_{\sigma(t)}(t) \right| \\
        & \leq & \left|x(t)-\hat{x}_{\bar{\mathcal{P}}} \right| + \left| \hat{x}_{\bar{\mathcal{P}}} - \hat{x}_{\sigma(t)}(t) \right| \\
         & \leq & \left|x(t)-\hat{x}_{\bar{\mathcal{P}}} \right| + \pi_{\sigma(t)}(t) \\
         & \leq & \left|x(t)-\hat{x}_{\bar{\mathcal{P}}} \right| + \pi_{\bar{\mathcal{I}}}(t),
    \end{array}
\end{equation}
where we have obtained the second to last inequality with \eqref{eq:pi_sigma} and the last inequality by definition of $\pi_{\sigma(t)}$.

Therefore, using \eqref{eq:pi_Ibar},  \eqref{eq:tildex_X} and \eqref{eq:tildex_Pbar}, 
\begin{equation}
    \begin{array}{ll}
\left|x(t)-\hat{x}_{\sigma(t)}(t)\right|  \leq & \beta_{\bar{\mathcal{P}}}\left(\left| \hat{x}_{\bar{\mathcal{P}}}(0)-x(0) \right|,t \right) \\
          & + \beta_{\bar{\mathcal{I}}}\left(\left| \hat{x}_{\bar{\mathcal{I}}}(0)-x(0) \right|,t \right) \\
          & + \bar{\beta}_{{\mathcal{P}}}\left(\left| \hat{x}_{{\mathcal{P}}}(0)-x(0) \right|,t \right).
    \end{array}
\end{equation}
where $\bar{\beta}:=\underset{\mathcal{P}\subset     \bar{\mathcal{I}}:\#(\mathcal{P})=N-2M}{\max} \beta_{\mathcal{P}}$ which is also a $\mathcal{KL}$ function.

By a property of class $\mathcal{KL}$ functions \cite[Lemma 5.3]{kellett2004weak}, we know that there exist class $\mathcal{K}_{\infty}$ functions $\alpha_{\bar{\mathcal{P}}}$, $\alpha_{\bar{\mathcal{I}}}$ and $\alpha_{{\mathcal{P}}}$, as well as positive constants $\lambda_{\bar{\mathcal{P}}}$, $\lambda_{\bar{\mathcal{I}}}$ and 
$\lambda_{{\mathcal{P}}}>0$ such that
\begin{equation}
    \begin{array}{ll}
\left|x(t)-\hat{x}_{\sigma(t)}(t)\right|  \leq & \alpha_{\bar{\mathcal{P}}}\left(\left| \hat{x}_{\bar{\mathcal{P}}}(0)-x(0) \right|\right) e^{-\lambda_{\bar{\mathcal{P}}} t} \\
          & + \alpha_{\bar{\mathcal{I}}}\left(\left| \hat{x}_{\bar{\mathcal{I}}}(0)-x(0) \right|\right) e^{-\lambda_{\bar{\mathcal{I}}} t} \\
          & + \alpha_{{\mathcal{P}}}\left(\left| \hat{x}_{{\mathcal{P}}}(0)-x(0) \right|\right) e^{-\lambda_{{\mathcal{P}}} t}, \, \forall t\geq0.
    \end{array}
\end{equation}

Hence, we obtain \eqref{eq:tildex_algo} as desired with $\tilde{\beta}:=3 \max \left\{ \alpha_{\bar{\mathcal{P}}}\left(\left| \hat{x}_{\bar{\mathcal{P}}}(0)-x(0) \right|\right),\, \alpha_{\bar{\mathcal{I}}}\left( \left| \hat{x}_{\bar{\mathcal{I}}}(0)-x(0) \right| \right),\, \right.$ $ \left.\alpha_{{\mathcal{P}}}\left(\left| \hat{x}_{{\mathcal{P}}}(0)-x(0) \right|\right)  \right\} e^{-\lambda t} $, where $\lambda := \min \left\{ \lambda_{\bar{\mathcal{P}}}, \lambda_{\bar{\mathcal{I}}}, \lambda_{{\mathcal{P}}} \right\}$.

\subsection{Proof of Proposition \ref{prop:spec}}
Let the state estimation error for every set $\mathcal{J}\subset\mathbb{N}_{[1,N]}$ be $
\tilde{x}_{\mathcal{J}}:=x-\hat{x}_{\mathcal{J}}$. The state estimation error system is
\begin{align}
    \dot{\tilde{x}}_{\mathcal{J}} & = \left(A - L_{\mathcal{J}} H_{\mathcal{J}} \right) \tilde{x}_{\mathcal{J}} + \phi(z) - \phi(\xi_{\mathcal{J}}) - L_{\mathcal{J}} a_{\mathcal{J}}.
\end{align}
Since the nonlinearity $\phi$ satisfies Assumption \ref{assum:sector}, there exists $\delta_i(t)\in\left[\underline{d}_{i},\bar{d}_{i}\right]$, for $i\in\mathbb{N}_{[1,N]}$ such that $$\phi(z)-\phi(\xi_{\mathcal{J}})=\delta(t)\eta_{\mathcal{J}}-\delta(t)K_{\mathcal{J}}a_{\mathcal{J}},$$ 
where $\delta(t):=\textrm{diag}\left(\delta_1(t),\delta_2(t),\dots,\delta_{N}(t) \right)$ and $\eta_{\mathcal{J}}:=\left(H-K_{\mathcal{J}}H_{\mathcal{J}} \right)\tilde{x}_{\mathcal{J}}$. Then the state estimation error system is
\begin{equation}
    \dot{\tilde{x}}_{\mathcal{J}}  = \left(A - L_{\mathcal{J}} H_{\mathcal{J}} \right) \tilde{x}_{\mathcal{J}} + 
    \delta(t) \eta_{\mathcal{J}} - \left( \delta(t) K_{\mathcal{J}} + L_{\mathcal{J}} \right) a_{\mathcal{J}}. \label{eq:state_error_sys}
\end{equation}

We are now ready to show that observer \eqref{eq:xhat_J} satisfies \eqref{eq:obs_iss} of Theorem \ref{thm:M_obs} using a candidate Lyapunov function $V_{\mathcal{J}}(\tilde{x}_{\mathcal{J}})=\tilde{x}_{\mathcal{J}}^{T}P_{\mathcal{J}}\tilde{x}_{\mathcal{J}}$, where $P_{\mathcal{J}}=P_{\mathcal{J}}^{T}>0$ satisfies \eqref{eq:cc_lmi}. The time derivative of $V_{\mathcal{J}}(\tilde{x}_{\mathcal{J}})$ along the solutions of the state estimation error system \eqref{eq:state_error_sys} is
\begin{equation}
     \dot{V}  \left(\tilde{x}_{\mathcal{J}}\right)  =  \chi_{\mathcal{J}}^{T}  \left[\begin{array}{ccc} P_{\mathcal{J}} \tilde{A}_{\mathcal{J}} + \tilde{A}_{\mathcal{J}}^{T} P_{\mathcal{J}} & P_{\mathcal{J}} & -P_{\mathcal{J}} \\
    P_{\mathcal{J}} & 0 & 0 \\
    -P_{\mathcal{J}} & 0 & 0\end{array}\right] \chi_{\mathcal{J}},
\end{equation}
where $\chi_{\mathcal{J}} = \left( \tilde{x}_{\mathcal{J}}, \delta(t) \eta_{\mathcal{J}}, \left(\delta(t)K_{\mathcal{J}} + L_{\mathcal{J}} \right)a_{\mathcal{J}}\right)$ and $\tilde{A}_{\mathcal{J}}:=A - L_{\mathcal{J}} H_{\mathcal{J}} $.

Applying \eqref{eq:cc_lmi}, we obtain
\begin{align}
     \dot{V}  \left(\tilde{x}_{\mathcal{J}}\right) \leq &  -\nu_{\mathcal{J}}|\tilde{x}_{\mathcal{J}}|^2 - 2 \eta_{\mathcal{J}}^{T} \delta(t) \eta_{\mathcal{J}} \nonumber \\
    & + 2 \eta_{\mathcal{J}}^{T} \delta(t)^{T} \textrm{diag}\left(1/\bar{d}_{1}, 1/\bar{d}_{2}, \dots, 1/\bar{d}_{N} \right) \delta(t) \eta_{\mathcal{J}} \nonumber \\
    & + \mu_{a}\left|\delta(t)K_{\mathcal{J}} + L_{\mathcal{J}} \right|^2 \left|a_{\mathcal{J}} \right|^2.
\end{align}
We first examine the second and third terms on the right hand side of the inequality component by component, i.e. for $i\in\mathbb{N}_{[1,N]}$, $\delta_i(t)-\delta_i(t)^2/\bar{d}_{i}= \delta_i(t) \left(1-\delta_i(t)/\bar{d}_{i} \right) \geq 0$, as $\delta_i(t) > 0$ and $1-\delta_i(t)/\bar{d}_{i}\geq 0$ by Assumption \ref{assum:sector}. Next, $\left|\delta(t)K_{\mathcal{J}} + L_{\mathcal{J}} \right|^2  \leq 2 \left| \delta(t)K_{\mathcal{J}} \right|^2 + 2 \left| L_{\mathcal{J}} \right|^2$ by Young's inequality and we obtain 
\begin{align}
     \dot{V}  \left(\tilde{x}_{\mathcal{J}}\right) \leq &  -\nu_{\mathcal{J}}|\tilde{x}_{\mathcal{J}}|^2 + \mu_{a} \left( 2 \bar{d}^2 \left| K_{\mathcal{J}} \right|^{2} +   2 \left| L_{\mathcal{J}} \right|^{2} \right) \left|a_{\mathcal{J}} \right|^2, \label{eq:vbar_int}
\end{align}
where $\bar{d}:=\max\left\{\bar{d}_{1},\bar{d}_{2},\dots,\bar{d}_{N}\right\}$.

By noting that $V_{\mathcal{J}}(\tilde{x}_{\mathcal{J}})$ can be sandwiched as follows
\begin{equation}
    \lambda_{\min}\left( P_{\mathcal{J}} \right) \left|\tilde{x}_{\mathcal{J}} \right|^{2} \leq V_{\mathcal{J}}(\tilde{x}_{\mathcal{J}}) \leq \lambda_{\max}\left( P_{\mathcal{J}} \right) \left|\tilde{x}_{\mathcal{J}} \right|^{2} , \label{eq:sand_V}
\end{equation}
and using the comparison principle, the solution to \eqref{eq:vbar_int} satisfies
\begin{align}
    V\left(\tilde{x}_{\mathcal{J}}\right) & \leq e^{-\lambda_{\mathcal{J}}t} V \left( \tilde{x}_{\mathcal{J}}(0) \right) + {\alpha_{\mathcal{J}}} \int_{0}^{t} e^{-\lambda_{\mathcal{J}}(t-s)} \left| a_{\mathcal{J}}(s) \right|^2 ds,
\end{align}
where $\lambda_{\mathcal{J}}:=\nu_{\mathcal{J}}/\lambda_{\max}(P_{\mathcal{J}})$, $\alpha_{\mathcal{J}}:=2\mu_{a}\left(\bar{d}^2\left| K_{\mathcal{J}}\right|^2 + \left| L_{\mathcal{J}}\right|^2\right)/\lambda_{\max}\left(P_{\mathcal{J}} \right)$. Since $\int_{0}^{t} e^{-\lambda_{\mathcal{J}}(t-s)}  ds = \left(1 - e^{-\lambda_{\mathcal{J}}t}\right)/\lambda_{\mathcal{J}} \leq 1/\lambda_{\mathcal{J}}$, we obtain
\begin{align}
    V\left(\tilde{x}_{\mathcal{J}}\right) & \leq e^{-\lambda_{\mathcal{J}}t} V \left( \tilde{x}_{\mathcal{J}}(0) \right) + {\alpha_{\mathcal{J}}/\lambda_{\mathcal{J}}} \underset{s\in[0,t]}{\sup} \left| a_{\mathcal{J}}(s) \right|^2 .
\end{align}
By applying \eqref{eq:sand_V} again, we obtain \eqref{eq:obs_iss} as desired with $\hat{\beta}(r,t)=\sqrt{\frac{\lambda_{\max}\left(P_{\mathcal{J}}\right)}{\lambda_{\min}\left(P_{\mathcal{J}}\right)}} e^{-\frac{\lambda_{\mathcal{J}}}{2} t} . r$ and $\hat{\gamma}(r):= \sqrt{\frac{\alpha_{\mathcal{J}}}{\lambda_{\mathcal{J}}\lambda_{\min}\left(P_{\mathcal{J}} \right)}} r$.

\subsection{Proof of Proposition \ref{prop:voltage}}
From \eqref{eq:line_rel}, \eqref{eq:droop} and \eqref{eq:vhat}, we have for $i\in\mathbb{N}_{[1,N]}$
\begin{equation}
    \begin{array}{lll}
    \left|v_{i}^{2}(t)-\hat{v}_{i}^{2}(t)\right| & \leq & \left| H_i \right| \left|{x}(t) - \hat{x}(t) \right|, \\
    &\leq & \left| H_i \right| \tilde{\beta}\left( \left|{x}(0) - \hat{x}(0) \right|, t \right),\; t\geq 0, 
    \end{array}
\end{equation}
where we obtain the last inequality according to Proposition \ref{prop:spec} and Theorem \ref{thm:algo}, and $\tilde{\beta}$ is a class $\mathcal{KL}$ function from Theorem \ref{thm:algo}. Hence, we have shown \eqref{eq:v2_vhat2} with $\beta_v(r,s) = \left|H_i \right| \tilde{\beta}(r,s)$.

\bibliographystyle{ieeetr}
\bibliography{grid_model_analysis.bib}

\end{document}